\documentclass[11pt]{amsart}
\usepackage{amsfonts, amssymb, amsmath, amsthm, hyperref, color, float,enumerate}
\usepackage{tikz,pgffor,amsthm}
\theoremstyle{plain}
   \newtheorem{teo}{Theorem}
   \newtheorem{coro}[teo]{Corollary}
   \newtheorem{lema}[teo]{Lemma}
   \newtheorem{propo}[teo]{Proposition}
   
\theoremstyle{definition}
   
\theoremstyle{remark}
 \newtheorem{obs}{Remark}

\numberwithin{equation}{section}
\allowdisplaybreaks
\newcommand{\R}{\mathbb{R}} 
\newcommand{\supp}{\mathrm{supp}}

\begin{document}

\title[Mixed inequalities involving radial functions]{Better bounds on mixed inequalities involving radial functions and applications}

\author[F. Berra]{Fabio Berra}
\address{CONICET and Departamento de Matem\'{a}tica (FIQ-UNL),  Santa Fe, Argentina.}
\email{fberra@santafe-conicet.gov.ar}

\thanks{The author was supported by CONICET and UNL}

\subjclass[2010]{42B25, 26A33}

\keywords{Maximal operators, Young functions, fractional integral operator}

\begin{abstract}
We prove mixed inequalities for the generalized maximal operator $M_\Phi$ when the function $v$ is a radial power function that fails to be locally integrable. Concretely, let $u$ be a weight, $v(x)=|x|^\beta$ with $\beta<-n$ and $r\geq 1$. If $\Phi$ is a Young function with certain properties, then the inequality
\[uv^r\left(\left\{x\in\mathbb{R}^n: \frac{M_\Phi (fv)(x)}{v(x)}>t\right\}\right)\leq C\int_{\mathbb{R}^n}\Phi\left(\frac{|f(x)|}{t}\right)v^r(x)Mu(x)\,dx\]
holds for every $t>0$ and every bounded function. This improves a similar mixed estimate proved in \cite{BCP-M}. 

As an application, we give mixed estimates for the generalized fractional maximal operator $M_{\gamma,\Phi}$, where $0<\gamma<n$ and $\Phi$ is of $L\log L$ type. A special case involving the fractional maximal operator $M_\gamma$ allows to obtain a similar estimate for the fractional integral operator $I_\gamma$ through an extrapolation result. Furthermore, we also give mixed estimates for commutators of singular integral Calder\'on-Zygmund operators and of $I_\gamma$, both with Lipschitz symbol. 
\end{abstract}

\maketitle

\section{Introduction}

In 1972 Muckenhoupt characterized the weights $w$ for which the Hardy-Littlewood maximal operator $M$ is bounded in $L^p(w)$ when $1<p<\infty$. The classes of such weights were called $A_p$ and acquired a primal role in many estimates of Harmonic Analysis.

Later on, Sawyer proved in \cite{Sawyer} an inequality involving two weights on the real line. More precisely, it was shown that if $u$ and $v$ are $A_1$ weights then the estimate
\begin{equation}\label{eq: intro - mixta de Sawyer}
uv\left(\left\{x\in\mathbb{R}: \frac{M(fv)(x)}{v(x)}>t\right\}\right)\leq \frac{C}{t}\int_{\mathbb{R}}|f(x)|u(x)v(x)\,dx
\end{equation}

holds for every positive $t$. This expression can be seen as the weak $(1,1)$ type of the operator $Sf=M(fv)/v$ with respect to the measure $d\mu(x)=u(x)v(x)\,dx$, and when we take $v=1$ we get the weak $(1,1)$ type for $M$ with respect to $A_1$ weights. A motivation for considering \eqref{eq: intro - mixta de Sawyer} is the fact that it provides a different proof of Muckenhoupt theorem. Indeed, it is not difficult to see that $S$ is bounded in $L^\infty(uv)$. By combining this fact with \eqref{eq: intro - mixta de Sawyer}, Marcinkiewicz interpolation theorem implies that $S$ is of strong $(p,p)$ type with respect to $uv$ for every $1<p<\infty$. For a fixed $p$ and $w\in A_p$ we can use Jones factorization theorem to write $w=uv^{1-p}$ where $u$ and $v$ are $A_1$ weights. Then
\[\|Mf\|_{L^p(w)}^p=\|S(fv^{-1})\|_{L^p(uv)}^p\leq C_p\|fv^{-1}\|_{L^p(uv)}^p=C_p\|f\|_{L^p(w)}^p.\]

In fact, $M$ is bounded in $L^p(w)$ if and only if $S$ is bounded in $L^p(uv)$, where $w=uv^{1-p}$. 

Further extensions of \eqref{eq: intro - mixta de Sawyer} were obtained later for different authors. In \cite{CruzUribe-Martell-Perez} mixed inequalities for $M$ and Calderón-Zygmund operators (CZO) in higher dimensions were obtained. Some extensions for commutators of CZO and operators with kernels of H\"{o}rmander type were proved in \cite{BCP-M}. 

Regarding to more general operators, such as the maximal function $M_\Phi$ defined by means of the Luxemburg averages $\|\cdot\|_{\Phi,Q}$, a first mixed inequality was proved in \cite{Berra}. We give the precise statement in the following theorem.

\begin{teo}[\cite{Berra}, Theorem 1]\label{teo: intro - B19}
	Let $r\geq 1$, $\delta\geq 0$ and $\Phi(t)=t^r(1+\log^+t)^\delta$. If $u, v^r$ are $A_1$ weights and $w=1/\Phi(1/v)$, then the inequality
	\[uw\left(\left\{x\in\mathbb{R}^n: \frac{M_\Phi (fv)(x)}{v(x)}>t\right\}\right)\leq C\int_{\mathbb{R}^n}\Phi\left(\frac{|f(x)|v(x)}{t}\right)u(x)\,dx\]
	holds for every positive $t$.
\end{teo}

Although the theorem above provided a newfangled estimate for these type of operators, it had the disadvantage that it failed to be homogeneous in $v$. This problem was overcome in \cite{BCP-MN}, where we improved the estimate above and considered a bit wider family of Young functions than just $L\log L$.  Given $r\geq 1$, we say that a Young function belongs to the family $\mathfrak{F}_r$ if $\Phi$ is submultiplicative, has lower type $r$ and satisfies the condition
\begin{equation}\label{eq: intro - condicion Fr}
\frac{\Phi(t)}{t^r}\leq C_0 (\log t)^\delta
\end{equation}
for some constants $C_0>0$, $\delta\geq 0$ and every $t\geq t^*\geq 1$. The precise statement of the improvement of Theorem~\ref{teo: intro - B19} is given below.

\begin{teo}[\cite{BCP-MN}, Theorem 1.1]\label{teo: intro - mixta para MPhi BCP(MN)}
	Let $r\geq 1$, $\Phi\in \mathcal{F}_r$ and $u, v^r$ be $A_1$ weights. Then there exists a positive constant $C$ such that for every $t>0$ the inequality
	\[uv^r\left(\left\{x\in\mathbb{R}^n: \frac{M_\Phi (fv)(x)}{v(x)}>t\right\}\right)\leq C\int_{\mathbb{R}^n}\Phi\left(\frac{|f(x)|}{t}\right)u(x)v^r(x)\,dx\]
	holds.
\end{teo} 

Further extensions and generalizations of the results above can also be found in \cite{L-O-P} and \cite{Berra21}. 

\medskip

In this paper we study mixed inequalities where $v$ is a radial power function in $\mathbb{R}^n$ that is not even locally integrable. The fact that one can prove such estimates under this condition may be surprising. However, many inequalities of this type have been already proved in the literature. In \cite{Andersen-Muckenhoupt} the authors showed that if $u\in A_1$ and $v(x)=|x|^{-d}$ with $d\neq 1$ then the inequality
\[uv\left(\left\{x\in\mathbb{R}: \frac{|\mathcal{T}f(x)|}{v(x)}>t\right\}\right)\leq \frac{C_d}{t}\int_{\mathbb{R}}|f(x)|u(x)\,dx\] 
holds for every $t>0$, where $\mathcal{T}$ is either $M$ or the Hilbert transform. 
A similar version of this estimate for CZO in higher dimensions was obtained in \cite{MR-OS-SG}. 

Later on, Ombrosi and Pérez proved in \cite{O-P} mixed inequalities for the Hardy-Littlewood maximal operator in $\mathbb{R}^n$. More precisely, they showed that if $u$ is any nonnegative locally integrable function and $v(x)=|x|^{-nr}$, with $r>1$ then the inequality
\[uv\left(\left\{x\in\mathbb{R}^n: \frac{M(fv)(x)}{v(x)}>t\right\}\right)\leq \frac{C}{t}\int_{\mathbb{R}^n}|f(x)|Mu(x)v(x)\,dx\]
holds. An extension of this result for $M_\Phi$ when $\Phi$ is a function of $L\log L$ type was obtained in \cite{BCP-M}, as we state in the following theorem. 
\begin{teo}[\cite{BCP-M}, Theorem 5]\label{teo: mixta para Mphi BCP(M)}
	Let $\Phi(t)=t^r(1+\log^+t)^\delta$, where $r\geq 1$ and $\delta\geq 0$. Let $u$ be any weight, $v(x)=|x|^\beta$ with $\beta<-n$ and $w=1/\Phi(1/v)$.  Then there exists a positive constant $C$ such that the inequality
\[uw\left(\left\{x\in\mathbb{R}^n: \frac{M_\Phi(fv)(x)}{v(x)}>t\right\}\right)\leq C\int_{\mathbb{R}^n}\Phi\left(\frac{|f(x)|v(x)}{t}\right)Mu(x)\,dx\]	
holds for every $t>0$.
\end{teo}

We can see that this last estimate is similar to the corresponding one in Theorem~\ref{teo: intro - B19}. Thus it also fails to be homogeneous in $v$ and, in views of the improvement given by Theorem~\ref{teo: intro - mixta para MPhi BCP(MN)}, it is natural to think whether a similar estimate can be achieved when $v$ is a radial power function. In this article we shall answer this question positively: by adapting some of the ideas involved in both Theorem~\ref{teo: intro - mixta para MPhi BCP(MN)} and~\ref{teo: mixta para Mphi BCP(M)} we can prove the following, which is our main result.
 
 \begin{teo}\label{teo: desigualdad mixta para MPhi}
 	Let $u$ be any weight and $v(x)=|x|^{\beta}$, with $\beta<-n$. Let $r\geq 1$ and $\Phi$ be a Young function belonging to $\mathfrak{F}_r$. Then there exists a positive constant $C$ such that the inequality
 	\[uv^r\left(\left\{x\in \mathbb{R}^n: \frac{M_\Phi(fv)(x)}{v(x)}>t\right\}\right)\leq C\int_{\mathbb{R}^n}\Phi\left(\frac{|f(x)|}{t}\right)v^r(x)Mu(x)\,dx\] 
 	holds for every $t>0$ and every bounded function.
 \end{teo}

As we said above, this last estimate improves Theorem~\ref{teo: mixta para Mphi BCP(M)} since we are getting a mixed inequality that is homogeneous in $v$. Observe that every Young function of the type $\Phi(t)=t^r(1+\log^+t)^\delta$ with $r\geq 1$ and $\delta\geq 0$ belongs to $\mathfrak{F}_r$. Therefore we have the following corollary.

\begin{coro}\label{coro: mixta para L log L}
	Let $u$ be any weight and $v(x)=|x|^{\beta}$, with $\beta<-n$. If $\Phi(t)=t^r(1+\log^+t)^\delta$ where $r\geq 1$ and $\delta\geq 0$, then there exists a positive constant $C$ such that the inequality
	\[uv^r\left(\left\{x\in \mathbb{R}^n: \frac{M_\Phi(fv)(x)}{v(x)}>t\right\}\right)\leq C\int_{\mathbb{R}^n}\Phi\left(\frac{|f(x)|}{t}\right)v^r(x)Mu(x)\,dx\] 
	holds for every $t>0$ and every bounded function.
\end{coro}

If $w=1/\Phi(1/v)$, we can easily check that $uw\leq v^r$ and $\Phi(|f|/t)v^r\leq \Phi(|f|v/t)$. So we get better bounds on this last estimate. 

\medskip

\begin{obs}
As we have mentioned, mixed inequalities can be applied to obtain strong types of the operator involved on weighted $L^p$ spaces. For instance, in \cite{BCP-MN} we proved that Theorem~\ref{teo: intro - mixta para MPhi BCP(MN)} implies that the operator $S_\Phi f=M_\Phi(fv)/v$ is bounded on $L^p(uv^r)$ (and it is, consequently, of weak $(p,p)$ type), for $r<p<\infty$. This allowed us to give a different proof of the fact that $M_\Phi: L^p(w)\hookrightarrow L^p(w)$ for every $r<p<\infty$ and $w\in A_{p/r}$, when the function $\Phi$ belongs to the class $\mathfrak{F}_r$. The techniques encompass interpolation for modular type inequalities and the fact that $M_\Phi$ is bounded on $L^\infty(uv^r)$. Observe that this last property fails when we consider the radial power function $v(x)=|x|^\beta$. This discussion leads to the question if we can prove other weak or strong $(p,p)$ types for $S_\Phi$ in this setting. The answer is negative in general. To see this, assume for simplicity that $\Phi(t)=t$ and that there exists $1<p_0<\infty$ such that $S_\Phi$ is of weak $(p_0,p_0)$ type. Then, by Marcinkiewicz interpolation theorem, we would have that $S_\Phi$ is of strong $(p,p)$ type with respect to $uv$, for every $1<p<p_0$. This would yield that $M$ is bounded on $L^p(w)$ for $w=uv^{1-p}$. By Muckenhoupt theorem we would get that $uv^{1-p}\in A_p$. If we take $u=1$ this leads to $|x|^{\beta(1-p)}\in A_p$, a contradiction. This shows that weak $(p,p)$ types might not hold when $v$ is the considered radial power function. Thus we shall give endpoint mixed type inequalities only.  	
\end{obs}

As a consequence of Corollary~\ref{coro: mixta para L log L} we can prove mixed inequalities involving radial functions for  the \emph{generalized fractional maximal operator}, defined by
\[M_{\gamma,\varphi}f(x)=\sup_{Q\ni x} |Q|^{\gamma/n}\|f\|_{\varphi, Q},\]
where $0<\gamma<n$ and $\varphi$ is a Young function. By using a pointwise estimate of Hedberg type (see Proposition~\ref{propo: estimacion puntual M_{gamma, Phi}} for details) we can obtain the following result.

\begin{teo}\label{teo: mixta para Mgamma,phi}
	Let $0<\gamma<n$, $1\leq r<n/\gamma$ and $q$ such that $1/q=1/r-\gamma/n$. For $\delta\geq 0$ define $\Phi(t)=t^r(1+\log^+t)^\delta$. If $u$ is a weight and $v(x)=|x|^\beta$ with $\beta<-n$, there exists a positive constant $C$ such that the inequality
	\[uv^q\left(\left\{x\in\mathbb{R}^n: \frac{M_{\gamma,\Phi}(fv)(x)}{v(x)}>t\right\}\right)\leq C\varphi\left(\int_{\mathbb{R}^n} \Phi_\gamma\left(\frac{|f(x)|}{t}\right)\Psi((Mu)^{1/q}(x)v(x))\,dx\right)\]
	holds for every bounded function, where $\varphi(t)=t^{q/r}(1+\log^+t)^{\delta q/r}$, $\Psi(t)=t^r(1+\log^+(t^{1-q/r}))^{n\delta/(n-r\gamma)}$ and $\Phi_\gamma(t)=\Phi(t)(1+\log^+t)^{\delta r\gamma/(n-r\gamma)}$.
\end{teo}

When $u$ and $v^q$ are $A_1$ weights, this inequality was obtained in \cite{FB}. A further extension when $v^q\in A_\infty$ appears in \cite{Berra21}.

Theorem~\ref{teo: mixta para Mgamma,phi} gives a quite general mixed inequality that implies some interesting particular cases. If we take $\Phi(t)=t$, that is, we pick $r=1$ and $\delta=0$ we have an estimate for the classical fractional maximal operator $M_\gamma$. For the next results we shall assume that the weight $u$ belongs to $A_1$.

\begin{coro}\label{coro: mixta para M_gamma}
	Let $0<\gamma<n$, $u\in A_1$, $v(x)=|x|^\beta$ with $\beta<-n$ and $q=n/(n-\gamma)$. Then there exists a positive constant $C$ such that the inequality
	\[uv^q\left(\left\{x\in\mathbb{R}^n: \frac{M_\gamma(fv)(x)}{v(x)}>t\right\}\right)^{1/q}\leq \frac{C}{t}\int_{\mathbb{R}^n}|f(x)|u^{1/q}(x)v(x)\,dx\]
	holds for every bounded function. 
\end{coro}

When $\gamma=0$ we have $q=1$ and $M_\gamma=M$. This case corresponds to the mixed estimate obtained in \cite{O-P}.

The corollary above allows us to obtain a similar estimate for the fractional integral operator $I_\gamma$, defined by
\[I_\gamma f(x)=\int_{\mathbb{R}^n} \frac{f(y)}{|x-y|^{n-\gamma}},\]
for every $f$ such that the integral is finite in almost every $x$. There is a well-known relation between $I_\gamma$ and $M_\gamma$ proved by Muckenhoupt and Wheeden in \cite{Muckenhoupt-Wheeden74}: for every $0<p<\infty$ and every $w\in A_\infty$ the inequality
\begin{equation}\label{eq: intro - tipo Coifman Igamma-Mgamma}
\int_{\mathbb{R}^n}|I_\gamma f(x)|^pw(x)\,dx\leq C\int_{\mathbb{R}^n}(M_\gamma f(x))^pw(x)\,dx
\end{equation}
holds. 

By adapting an extrapolation argument from \cite{O-P} (see Theorem~\ref{teo: extrapolacion}) and combining it with \eqref{eq: intro - tipo Coifman Igamma-Mgamma} we can prove the following theorem.

\begin{teo}\label{teo: mixta para I_gamma}
	Let $0<\gamma<n$, $u\in A_1$, $v(x)=|x|^\beta$ with $\beta<-n$ and $q=n/(n-\gamma)$. Then there exists a positive constant $C$ such that the inequality
	\[uv^q\left(\left\{x\in\mathbb{R}^n: \frac{|I_\gamma(fv)(x)|}{v(x)}>t\right\}\right)^{1/q}\leq \frac{C}{t}\int_{\mathbb{R}^n}|f(x)|u^{1/q}(x)v(x)\,dx\]
	holds for every bounded function. 
\end{teo}

The article is organized as follows. In Section~\ref{section: preliminares} we give the definitions and basic results that we will need for the main proofs, which are contained in Section~\ref{section: demostraciones}. We devote Section~\ref{seccion: caso Phi(t)=t} to prove Theorem~\ref{teo: mixta para I_gamma} by using an extrapolation technique. Finally, as an application of this last result, we obtain mixed inequalities for commutators with Lipschitz symbol in Section~\ref{seccion: aplicaciones}.

\section{Definitions and auxiliary lemmas}\label{section: preliminares}

Recall that a \emph{weight} $w$ is a locally integrable function such that $0<w(x)<\infty$ for almost every $x$. Given $1<p<\infty$, the \emph{Muckenhoupt $A_p$ class} is the collection of weights $w$
satisfying that there exists a positive constant $C$ such that the inequality
\[\left(\frac{1}{|Q|}\int_Q w\right)\left(\frac{1}{|Q|}\int_Q w^{-\frac{1}{p-1}}\right)^{p-1}
\leq C\] 
holds for every cube $Q\subseteq \R^n$  with sides parallel
to the coordinate axes.  When $p=1$ we say that $w\in A_1$ if there
exists a positive constant $C$ such that
\begin{equation*}
\frac{1}{|Q|}\int_Q w\le C\, \inf_Q w,
\end{equation*}
for every cube $Q\subseteq \R^n$. The $A_\infty$  class is defined as the collection of all the $A_p$ classes, that is, $A_\infty=\bigcup_{p\geq 1}A_p$. It is well-known that $A_p\subseteq A_q$ when  $p<q$. For further properties and the basic theory of Muckenhoupt weights we refer the reader to~\cite{javi}, ~\cite{GC-RdF} and~\cite{grafakos}.

By a \emph{Young function} we shall understand a mapping $\varphi:[0,\infty)\to[0,\infty]$ which is convex, strictly increasing, $\varphi(0)=0$ and verifies that $\varphi(t)\to\infty$ when $t\to\infty$. 

We say $\varphi$ has a \textit{lower type} $p$, $0<p<\infty$ if there exists a positive constant $C_p$ such that
\[\varphi(st)\leq C_ps^p\varphi(t),\]
for every $0<s\leq 1$ and $t>0$. As an immediate consequence of this definition we have that if $\varphi$ has lower type $p$, then $\varphi$ has lower type $\tilde p$ for every $0<\tilde p<p$. 

With $\varphi^{-1}$ we will denote the \emph{generalized inverse} of $\varphi$ given by the expression
\[\varphi^{-1}(t)=\inf\{s\geq 0: \varphi(s)\geq t\},\]
where we understand $\inf\emptyset=\infty$.

It is well-known that if $\varphi(t)=t^r(1+\log^+t)^\delta$, with $r\geq 1$ and $\delta\geq 0$ then
\begin{equation}\label{eq: preliminares - inversa de L log L}
\varphi^{-1}(t)\approx \frac{t^{1/r}}{(1+\log^+t)^{\delta/r}},
\end{equation}
where the notation $A\approx B$ means that there exist two positive constants $C_1$ and $C_2$ such that
\[C_1 A\leq B\leq C_2 A.\]

The \emph{generalized maximal operator} $M_\varphi$ is defined by
\begin{equation}\label{eq: definicion de MPhi}
M_\varphi f(x)=\sup_{Q\ni x} \|f\|_{\varphi, Q}
\end{equation}
where  $\|f\|_{\varphi, Q}$ denotes the Luxemburg type average of $f$ over $Q$ given by
\[\|f\|_{\varphi, Q}=\inf\left\{\lambda>0: \frac{1}{|Q|}\int_Q \varphi\left(\frac{|f(x)|}{\lambda}\right)\,dx\leq 1\right\}.\]

The following lemma establishes a control of the level sets of the operators above by dyadic cubes. A proof can be found in \cite{CUMP11}. 

\begin{lema}\label{lema: descomposicion de CZ del espacio}
	Given $\lambda>0$, a Young function $\varphi$ and a bounded function $f$, there exists a family of maximal dyadic cubes $\{Q_j\}$ that satisfies
	\[\{x\in\mathbb{R}^n: M_{\varphi}f(x)>\lambda\}\subseteq\bigcup_j (3Q_j),\]
	and $\left\|f\right\|_{\varphi,Q_j}>\lambda$  for every $j$.
\end{lema}

We finish this section with two technical results that will be useful in our main proof. The latter was set and proved in \cite{O-P}.

\begin{lema}\label{lema: estimacion constante epsilon}
	Let $f$ be the function defined in $[0,\infty)$ by 
	\[f(x)=\left\{\begin{array}{ccl}
	\left(1+\frac{1}{x}\right)^{\frac{x}{1+x}}& \textrm{ if } & x>0,\\
	1 & \textrm{ if} & x=0.
	\end{array}
	\right.\]
	Then we have that $1\leq f(x)\leq e^{1/e}$, for every $x\geq 0$.
\end{lema}

\begin{lema}\label{lema: igualdad de integral}
	Let $f$ be a positive and locally integrable function. Then for each $\ell,\lambda>0$ there exists a number $a\in\R^+$ which depends on $f$ and $\lambda$ that satisfies
	\bigskip
	\[\left(\int_{|y|\leq a^{\ell}}f(y)\,dy\right)a^n=\lambda.\]
\end{lema}

\section{Proofs of the main results}\label{section: demostraciones}
We devote this section to the proof of Theorem~\ref{teo: desigualdad mixta para MPhi} and Theorem~\ref{teo: mixta para Mgamma,phi}.
\begin{proof}[Proof of Theorem~\ref{teo: desigualdad mixta para MPhi}]
	Suppose without loss of generality that the right-hand side of the inequality is finite. This implies that $\Phi(|f|)v^r\in L^1$, which will yield that some quantities are finite. Indeed, we have that
	\[\int_{\mathbb{R}^n}\Phi(|f(x)|)v^r(x)\,dx=\int_{B(0,1)}+\int_{B(0,1)^c}.\]
	For the second term we can simply use polar coordinates to get
	\[\int_{B(0,1)^c}\Phi(|f(x)|)v^r(x)\,dx\leq C\int_1^\infty \rho^{r\beta+n-1}\,d\rho<\infty\]
	since $r\beta+n<0$. On the other hand, if $\tilde Q$ is any cube that contains $B(0,1)$ we can write
	\[\frac{u(\tilde Q)}{|\tilde Q|}\int_{B(0,1)}\Phi(|f(x)|)v^r(x)\,dx\leq \int_{\mathbb{R}^n}\Phi(|f(x)|)v^r(x)Mu(x)\,dx<\infty.\]
	
	It will also be enough to assume that $t=1$ by homogeneity. For $k\in\mathbb{Z}$ we define the sets
	\[G_k=\{x\in \mathbb{R}^n: 2^k<|x|\leq 2^{k+1}\},\quad I_k=\{x\in \mathbb{R}^n: 2^{k-1}<|x|\leq 2^{k+2}\},\]
	\[L_k=\{x\in \mathbb{R}^n: 2^{k+2}<|x|\}, \quad\textrm{ and }\quad C_k=\{x\in \mathbb{R}^n: |x|\leq 2^{k-1}\}.\]
	We also split $f=f\mathcal{X}_{I_k}+f\mathcal{X}_{L_k}+f\mathcal{X}_{C_k}=f_1+f_2+f_3$. It will be enough to achieve the estimate for each $f_i$, $1\leq i\leq 3$.
	
	We start with $f_1$. If we set $\Omega_k=\{x\in\mathbb{R}^n: M_{\Phi}(f_1v)(x)>2^{(k+1)\beta}\}$, then by Lemma~\ref{lema: descomposicion de CZ del espacio} we can find for every $k$ a disjoint collection of dyadic cubes $\{Q_j^k\}_j$ such that $\Omega_k\subseteq \bigcup_j 3Q_j^k$ and 
	\begin{equation}\label{eq: estimacion norma luxemburgo en Q_j}
	\|f_1v\|_{\Phi, Q_j^k}> 2^{(k+1)\beta} \quad\textrm{ for every }j.
	\end{equation}
	
	Then we proceed as follows
	\begin{align*}
	uv^r\left(\left\{x\in \mathbb{R}^n: M_{\Phi}(f_1v)(x)>v(x)\right\}\right)&\leq \sum_{k\in\mathbb{Z}}uv^r\left(\left\{x\in G_k: M_{\Phi}(f_1v)(x)>v(x)\right\}\right)\\
	&\leq  \sum_{k\in\mathbb{Z}} 2^{k\beta r}u\left(\left\{x\in G_k: M_{\Phi}(f_1v)(x)>2^{(k+1)\beta}\right\}\right)\\
	&=\sum_{k\in\mathbb{Z}} 2^{k\beta r}u\left(G_k\cap \Omega_k\right)\\
	&\lesssim \sum_{k\in\mathbb{Z}}\sum_j  2^{k\beta r}u(G_k\cap 3Q_j^k).
	\end{align*}

	For fixed $k$ and $j$, by \eqref{eq: estimacion norma luxemburgo en Q_j} we have that 
	\[\left\|\frac{f_1v}{2^{(k+1)\beta}}\right\|_{\Phi,Q_j^k}>1 \quad\textrm{ or equivalently}\quad 1<\frac{1}{|Q_j^k|}\int_{Q_j^k}\Phi\left(\frac{|f_1(x)|v(x)}{2^{(k+1)\beta}}\right)\,dx.\]
	We split now the cube $Q_j^k$ into the sets $A=\left\{x\in Q_j^k: v(x)\leq 2^{(k+1)\beta}t^*\right\}$ and $B=Q_j^k\backslash A$, where $t^*$ is the number involved in \eqref{eq: intro - condicion Fr}. Then
	\[1<\left\|\frac{f_1v}{2^{(k+1)\beta}}\right\|_{\Phi,Q_j}\leq \left\|\frac{f_1v\mathcal{X}_A}{2^{(k+1)\beta}}\right\|_{\Phi,Q_j}+\left\|\frac{f_1v\mathcal{X}_B}{2^{(k+1)\beta}}\right\|_{\Phi,Q_j}=I+II.\]
	The inequality above implies that either $I>1/2$ or $II>1/2$. If the first one holds, we can estimate
	\begin{align*}
	1&<\frac{1}{|Q_j^k|}\int_A \Phi\left(2\frac{|f_1|v}{2^{(k+1)\beta}}\right)\\
	&\leq \frac{\Phi(2t^*)}{(t^*)^r|Q_j^k|}\int_A \Phi(|f_1|)\left(\frac{v^r}{2^{(k+1)\beta r}}\right),
	\end{align*}
	since $\Phi$ is submultiplicative and has lower type $r$. Consequently,
	\[|Q_j^k|2^{(k+1)\beta r}\leq C\int_{Q_j^k} \Phi(|f_1(x)|)v^r(x)\,dx.\]
	If $II>1/2$, by using the decay condition of $\Phi$ we can write
	\begin{align*}
	1&<\frac{1}{|Q_j^k|}\int_B \Phi\left(\frac{2|f_1|v}{2^{(k+1)\beta}}\right)\\
	&\leq \frac{\Phi(2)C_0}{|Q_j^k|}\int_B \Phi(|f_1|)\left(\frac{v}{2^{(k+1)\beta}}\right)^r\left(\log\left(\frac{v}{2^{(k+1)\beta}}\right)\right)^\delta\\
	&\leq \frac{\Phi(2)C_0(-2\beta)\log 2}{|Q_j^k|}\int_{Q_j^k} \Phi(|f_1|)\left(\frac{v}{2^{(k+1)\beta}}\right)^r,
	\end{align*}
	since for $x\in I_k$ we have $v(x)\leq 2^{(k-1)\beta}$. This yields
	\[|Q_j^k|2^{(k+1)\beta r}\leq C\int_{Q_j^k} \Phi(|f_1(x)|)v^r(x)\,dx.\]
	By using these estimates, we have that
	\begin{align*}
	\sum_{k\in\mathbb{Z}}\sum_j 2^{k\beta r}u(3Q_j^k)&\leq C2^{-\beta r}\sum_{k\in\mathbb{Z}}\sum_j \frac{u(3Q_j^k)}{|3Q_j^k|}\int_{Q_j^k}\Phi(|f_1(x)|)v^r(x)\,dx\\
	&\leq C\sum_{k\in\mathbb{Z}} \int_{I_k}\Phi(|f(x)|)v^r(x)Mu(x)\,dx\\
	&\leq C\int_{\mathbb{R}^n}\Phi(|f(x)|)v^r(x)Mu(x)\,dx
	\end{align*}
	since the sets $I_k$ overlap at most three times.
	
	We turn now our attention to the estimate corresponding to $f_2$. Let $\theta$ be a positive number to be chosen later. We define the functions
	\[F_1(x)=\left(\int_{\{|y|>|x|\}}\frac{\Phi(|f_2(y)|)}{|y|^n}v^r(y)\,dy\right)^{1/r}\]
	and
	\[F_2(x)=\left[\left(\int_{\{|y|>|x|\}}\frac{\Phi(|f_2(y)|)}{|y|^n}v^r(y)\,dy\right)v^{\delta\theta}(x)\right]^{1/(r+\delta\theta)}.\]
	Set $\lambda_0(x)=\max\{F_1(x),F_2(x)\}$. We shall prove that
	\begin{equation}\label{eq: control de la norma de f_2v por lambda0}
	\|f_2v\|_{\Phi,B}\leq C\lambda_0(x).
	\end{equation}
	for every $x\in G_k$ and every ball $B$ that contains $x$. Fixed $x$ and $B$, if $y\in B\cap L_k$ we have $2|x|\leq 2^{k+2}<|y|$, so
	\[2R>|y-x|\geq |y|-|x|>\frac{|y|}{2},\]
	where $R$ is the radius of $B$.
	Let $\Gamma_1=\{z\in B: v(z)\leq \lambda_0 t^*\}$ and $\Gamma_2=B\backslash \Gamma_1$. Therefore
	\begin{align*}
	\frac{1}{|B|}\int_B \Phi\left(\frac{|f_2(y)|v(y)}{\lambda_0}\right)\,dy&\leq \frac{1}{|B|}\int_{\Gamma_1} \Phi\left(\frac{|f_2(y)|v(y)}{\lambda_0}\right)\,dy+\frac{1}{|B|}\int_{\Gamma_2} \Phi\left(\frac{|f_2(y)|v(y)}{\lambda_0}\right)\,dy\\
	&=I_1+I_2.
	\end{align*}
	Since $\Phi\in\mathfrak{F}_r$, we can estimate $I_1$ as follows
	\begin{align*}
	I_1&\leq \frac{\Phi(t^*)}{(\lambda_0t^*)^r|B|}\int_{\Gamma_1} \Phi(|f_2(y)|)v^r(y)\,dy\\
	&\leq \frac{C}{(F_1(x))^r}\int_{\{|y|>|x|\}} \frac{\Phi(|f(y)|)}{|y|^n}v^r(y)\,dy\\
	&\leq C.
	\end{align*}
	On the other hand, fix $\varepsilon>0$ and set $\gamma=1+\varepsilon$. By appyling Hölder inequality with $\gamma$ and $\gamma'$ and using the decay condition of $\Phi$ we get
	\begin{align*}
	I_2&\leq \frac{C_0}{|B|}\int_{\Gamma_2}\Phi(|f_2(y)|)\left(\frac{v(y)}{\lambda_0}\right)^r\log\left(\frac{v(y)}{\lambda_0}\right)^\delta\,dy\\
	&\leq \frac{C}{\lambda_0^r}\left(\frac{1}{|B|}\int_B\left[\Phi(|f_2(y)|)v^r(y)\right]^{\gamma}\,dy\right)^{1/\gamma}\left(\frac{1}{|B|}\int_{B\cap L_k} \left[\log\left(\frac{v}{\lambda_0}\right)\right]^{\delta\gamma'}\right)^{1/\gamma'}.
	\end{align*}
	Notice that the expression associated to $f_2$ is finite for every $\varepsilon$. Indeed, since $f$ is bounded we have
	\begin{align*}
	\int_B \left[\Phi(|f_2(y)|)v^r(y)\right]^{\gamma}\,dy&\leq \Phi(\|f\|_\infty)^\gamma\int_{|x|}^\infty\rho^{\beta r\gamma+n-1}\,d\rho\\
	&\leq \frac{\Phi(\|f\|_\infty)^\gamma}{-\beta r\gamma-n}|x|^{\beta r\gamma+n}.
	\end{align*}
	
	On the other hand, by using the fact that
	\begin{equation}\label{eq: teo: desigualdad mixta para MPhi - log t acotado por potencia}
	\log t\lesssim \xi^{-1}t^\xi \quad\textrm{ for every }\quad t,\xi>0
	\end{equation}
	we can estimate 
	\begin{align*}
	\left(\frac{1}{|B|}\int_{B\cap L_k} \left[\log\left(\frac{v}{\lambda_0}\right)\right]^{\delta\gamma'}\right)^{1/\gamma'}&\leq \frac{C}{\lambda_0^{\delta\theta}}\left(\frac{1}{\theta}\int_{\{|y|>|x|\}} |y|^{\beta\delta\theta\gamma'-n}\,dy \right)^{1/\gamma'}\\
	&\leq \frac{C}{\lambda_0^{\delta\theta}}\left(\frac{1}{\theta}\int_{|x|}^\infty \rho^{\beta\delta\theta\gamma'-1}\,d\rho \right)^{1/\gamma'}\\
	&\leq \frac{C}{\lambda_0^{\delta\theta}}\left(\frac{-1}{\theta^2\beta\delta\gamma'}\right)^{1/\gamma'}|x|^{\beta\theta\delta}\\
	&\leq \frac{C}{\lambda_0^{\delta\theta}}|x|^{\beta\delta\theta},
	\end{align*}
	provided we choose $\theta$ such that $\theta^2\beta\delta=-1$. Consequently, we have that
	\[I_2\leq \frac{C}{\lambda_0^{r+\delta\theta}}\left(\frac{1}{|B|}\int_B\left[\Phi(|f_2(y)|)v^r(y)\right]^{\gamma}\,dy\right)^{1/\gamma}v(x)^{\delta\theta}.\]
	From Lemma~\ref{lema: estimacion constante epsilon} we have that $C$ does not depend on $\varepsilon$. Thus, the expression above is true for every $\varepsilon>0$. By letting $\varepsilon\to 0$, we get that
	\begin{align*}
	I_2&\leq \frac{C}{\lambda_0^{r+\delta\theta}}\left(\int_{\{|y|>|x|\}}\frac{\Phi(|f_2(y)|)}{|y|^n}v^r(y)\,dy\right)v(x)^{\delta\theta}\\
	&\leq C
	\end{align*}
	because of the definition of $\lambda_0$. This allows us to conclude \eqref{eq: control de la norma de f_2v por lambda0}. If $G_k^1=G_k\cap \{x: \lambda_0(x)=F_1(x)\}$ and $G_k^2=G_k\cap \{x: \lambda_0(x)=F_2(x)\}$ we have that 
	
	\begin{align*}
	uv^r(\{x\in G_k: M_\Phi(f_2v)(x)>v(x)\})&\leq uv^r(\{x\in G_k^1: CF_1(x)>v(x)\})\\
	&\quad+uv^r(\{x\in G_k^2: CF_2(x)>v(x)\})
	\end{align*}
	and both of these terms can be bounded by
	\[C2^{k\beta r}u\left(\left\{x\in G_k: \int_{\{|y|>|x|\}}\frac{\Phi(|f(y)|)}{|y|^n}v^r(y)\,dy>v^r(x)\right\}\right).\]
	Therefore
	\begin{align*}
		uv^r\left(\left\{x: M_\Phi(f_2v)(x)>v(x)\right\}\right)&\leq C2^{k\beta r}u\left(\left\{x\in G_k: \int_{\{|y|>|x|\}}\frac{\Phi(|f(y)|)}{|y|^n}v^r(y)\,dy>v^r(x)\right\}\right)\\
	&\leq  C2^{k\beta r}u\left(\left\{x\in G_k: \int_{\{|y|>|x|\}}\frac{\Phi(|f(y)|)}{|y|^n}v^r(y)\,dy>2^{(k+1)\beta r}\right\}\right)\\
	&\leq C\int_0^\infty u\left(\left\{x\in G_k: \int_{\{|y|>|x|\}}\frac{\Phi(|f(y)|)}{|y|^n}v^r(y)\,dy>\tau\right\}\right)\,d\tau\\
	&=\int_{\mathbb{R}^n}\Phi(|f(y)|)v^r(y)\left(\frac{1}{|y|^n}\int_{\{|y|>|x|\}}u(x)\,dx\right)\,dy\\
	&\leq  \int_{\mathbb{R}^n}\Phi(|f(y)|)v^r(y)Mu(y)\,dy.
	\end{align*}
It only remains to show the desired bound for $f_3$. We shall see that
\begin{equation}\label{eq: control de la norma de f_3v}
\|f_3v\|_{\Phi, B}\leq C\left(\frac{1}{|x|^n}\int_{\{|y|\leq |x|/2\}}\Phi(|f_3(y)|)v^r(y)\,dy\right)^{1/r}=CF(x),
\end{equation}
for every $x\in G_k$ and every ball $B$ that contains $x$. Fix $x$, $B=B(x_B,R)$ and note that if $y\in C_k$ then
\[2R\geq |y-x|\geq |x|-|y|\geq \frac{|x|}{2}\geq |y|.\]
We split again the ball into the sets $\Gamma_1=\{z\in B: v(z)\leq t^*F(x)\}$ and $\Gamma_2=B\backslash \Gamma_1$. Thus,
\begin{align*}
\frac{1}{|B|}\int_B \Phi\left(\frac{|f_3(y)|v(y)}{F(x)}\right)\,dy&=\frac{1}{|B|}\int_{\Gamma_1}\Phi\left(\frac{|f_3(y)|v(y)}{F(x)}\right)\,dy+\frac{1}{|B|}\int_{\Gamma_2}\Phi\left(\frac{|f_3(y)|v(y)}{F(x)}\right)\,dy\\
&=I_1+I_2.
\end{align*}
For $I_1$ we use again the lower type of $\Phi$ to obtain
\[I_1\leq \frac{C|x|^{-n}}{(F(x))^r}\int_{\{|y|\leq |x|/2\}} \Phi(|f(y)|)v^r(y)\,dy\leq C.\]
On the other hand, since $\Phi\in\mathcal{F}_r$ we have that
\[I_2\leq C|x|^{-n}\int_{\{y\in\Gamma_2\,:\, |y|\leq |x|/2\}}\Phi(|f_3(y)|)\left(\frac{v(y)}{F(x)}\right)^r\log\left(\frac{v(y)}{F(x)}\right)^\delta\,dy=:\tilde{I}_2.\]
Set $A_j=\{y\in \Gamma_2: 2^{-j-1}|x|< |y|\leq 2^{-j}|x|\}$ for every $j\in\mathbb{N}$. Also define
\[\tilde{I}_2^N=C|x|^{-n}\int_{\bigcup_{j=1}^N A_j}\Phi(|f_3(y)|)\left(\frac{v(y)}{F(x)}\right)^r\log\left(\frac{v(y)}{F(x)}\right)^\delta\,dy\]
for every $N\in\mathbb{N}$. By fixing a positive integer $N$, $\varepsilon>0$ and $\gamma=1+\varepsilon$, Hölder inequality yields
\[\tilde{I}_2^N
\leq \frac{C|x|^{-n}}{(F(x))^r}\sum_{j=1}^N \left(\int_{A_j}\left[\Phi(|f_3(y)|)v^r(y)\right]^{\gamma}\,dy\right)^{1/\gamma}\left(\int_{A_j} \left[\log\left(\frac{v}{F(x)}\right)\right]^{\delta\gamma'}\right)^{1/\gamma'}\]

Notice that the integral associated to $f_3$ is finite for every $\varepsilon$ and every $j$. We fix $j$ and apply \eqref{eq: teo: desigualdad mixta para MPhi - log t acotado por potencia} with $\theta=-n/(2\beta\gamma'\delta)$ to estimate the last factor above as follows 
\begin{align*}
\left(\int_{A_j} \left[\log\left(\frac{v}{F(x)}\right)\right]^{\delta\gamma'}\right)^{1/\gamma'}&\leq \frac{C}{F(x)^{\delta\theta}}\left(\frac{1}{\theta}\int_{\{|y|\leq|x|/2\}} |y|^{\beta\delta\theta\gamma'}\,dy \right)^{1/\gamma'}\\
&\leq \frac{C}{F(x)^{\delta\theta}}\left(\frac{1}{\theta}\int_0^{|x|/2} \rho^{\beta\delta\theta\gamma'+n-1}\,d\rho \right)^{1/\gamma'}\\
&\leq \frac{C}{F(x)^{\delta\theta}}\left(\frac{1}{\theta}\frac{|x|^{n+\beta\delta\theta\gamma'}}{2^{n+\beta\delta\theta\gamma'}(n+\beta\delta\theta\gamma')}\right)^{1/\gamma'}\\
&\leq \frac{C}{F(x)^{\delta\theta}}\left(\frac{-4\beta\gamma'\delta}{2^{n/2}n^2}\right)^{1/\gamma'}|x|^{n/\gamma'+\beta\delta\theta}.
\end{align*}
Since $\gamma'\to\infty$ when $\varepsilon\to 0$, there exists $\varepsilon_0>0$ such that, if $0<\varepsilon<\varepsilon_0$ then 
\[\frac{-4\beta\delta}{2^{n/2}n^2}\leq \gamma'.\]
Consequently, by applying Lemma~\ref{lema: estimacion constante epsilon}, we get that for every $0<\varepsilon<\varepsilon_0$
\[\left(\int_{A_j} \left[\log\left(\frac{v}{F(x)}\right)\right]^{\delta\gamma'}\right)^{1/\gamma'}\leq \frac{C}{(F(x))^{\delta\theta}}e^{2/e}|x|^{n/(2\gamma')}.\]
By letting $\varepsilon\to 0$, we get
\[\left(\int_{A_j}\left[\Phi(|f_3(y)|)v^r(y)\right]^{\gamma}\,dy\right)^{1/\gamma}\left(\int_{A_j} \left[\log\left(\frac{v}{F(x)}\right)\right]^{\delta\gamma'}\right)^{1/\gamma'}\leq C\left(\int_{A_j} \Phi(|f(y)|)v^r(y)\,dy\right),\]
where $C$ does not depend on $\varepsilon$ nor $j$. This yields
\[\tilde{I}_2^N\leq \frac{C|x|^{-n}}{(F(x))^r}\left(\int_{\{|y|\leq |x|/2\}} \Phi(|f(y)|)v^r(y)\,dy\right)\leq C,\]
for every $N\in\mathbb{N}$ and with $C$ independent of $N$. The monotone convergence theorem finally implies \eqref{eq: control de la norma de f_3v}.

If we take $\ell=-n/(n+\beta r)>0$ and $\lambda=1$ in Lemma~\ref{lema: igualdad de integral}, there exists $a>0$ such that
\begin{equation}\label{eq: teo desigualdad mixta para MPhi - integral de Phi(f) igual a a^{-n}}
\int_{\{|y|\leq a^{\ell}\}}\Phi(|f(y)|)v^r(y)\,dy=a^{-n}.
\end{equation}
By combining these estimates we get
\begin{align*}
uv^r\left(\left\{x\in G_k: M_\Phi(f_3v)(x)>v(x)\right\}\right)&\leq uv^r\left(\left\{x\in G_k: CF(x)>v(x)\right\}\right)\\
&\leq uv^r\left(\left\{x: |x|\leq a^\ell, CF^r(x)>v^r(x)\right\}\right)+\\
&\quad +\sum_{j=0}^\infty uv^r\left(\left\{x: 2^ja^\ell<|x|\leq 2^{j+1}a^{\ell},\, CF^r(x)>v^r(x)\right\}\right)\\
&=A+B.
\end{align*}
By \eqref{eq: teo desigualdad mixta para MPhi - integral de Phi(f) igual a a^{-n}} we have that
\begin{align*}
A&=uv^r\left(\left\{x: |x|\leq a^\ell, \frac{C}{|x|^n}\int_{\{|y|\leq |x|/2\}}\Phi(|f(y)|)v^r(y)\,dy>v^r(x)\right\}\right)\\
&=uv^r\left(\left\{x: |x|\leq a^\ell, \frac{C}{|x|^n}a^{-n}>|x|^{r\beta}\right\}\right)\\
&=uv^r\left(\left\{x: ca^\ell<|x|\leq a^\ell\right\}\right),
\end{align*}
where $0<c<1$ since $C$ can be chosen greater than one. Therefore
\begin{align*}
A&\leq \int_{\{|x|>ca^\ell\}}u(x)v^r(x)\,dx\\
&=\sum_{j=1}^\infty \int_{\{2^{j-1}ca^\ell< |x|\leq 2^jca^\ell\}}u(x)v^r(x)\,dx\\
&\leq \sum_{j=1}^\infty (2^{j-1}ca^{\ell})^{r\beta}\int_{\{|x|< 2^jca^\ell\}}u(x)\,dx\\
&= C\sum_{j=1}^\infty (2^{j-1}ca^{\ell})^{r\beta+n}a^n\int_{\{|y|\leq a^\ell\}}\Phi(|f(y)|)v^r(y)\left(\frac{1}{(2^jca^\ell)^n}\int_{\{|x|< 2^jca^\ell\}}u(x)\,dx\right)\,dy\\
&\leq C\sum_{j=1}^\infty (2^{j-1}ca^{\ell})^{r\beta+n}a^n\int_{\{|y|\leq a^\ell\}}\Phi(|f(y)|)v^r(y)Mu(y)\,dy\\
&\leq C\int_{\mathbb{R}^n}\Phi(|f(y)|)v^r(y)Mu(y)\,dy,
\end{align*}
since $r\beta+n<0$ and therefore the series converges. Also observe that $C$ does not depend on $f$ since $\ell(r\beta+n)+n=0$. The estimate for $B$ is similar. Indeed, we have that
\begin{align*}
B&\leq\sum_{j=0}^\infty \int_{\{2^ja^\ell<|x|\leq 2^{j+1}a^\ell\}}u(x)v^r(x)\,dx\\
&\leq C\sum_{j=0}^\infty (2^ja^\ell)^{r\beta+n}a^n\int_{|y|\leq a^\ell}\Phi(|f(y)|)v^r(y)\left(\frac{1}{(2^{j+1}a^\ell)^n}\int_{\{|x|\leq 2^{j+1}a^\ell\}}u(x)\,dx\right)\,dy\\
&\leq C\int_{\mathbb{R}^n}\Phi(|f(y)|)v^r(y)Mu(y)\,dy.
\end{align*}
This completes the proof. \qedhere
\end{proof}

\medskip

Before turning to the proof of Theorem~\ref{teo: mixta para Mgamma,phi} we state a pointwise relation between $M_{\gamma,\Phi}$ and $M_\xi$ where $\Phi$ and $\xi$ are Young functions. This inequality can be seen as a Hedberg type estimate involving Luxemburg averages. A proof can be found in \cite{Berra21}.

\begin{propo}\label{propo: estimacion puntual M_{gamma, Phi}}
	Let $0<\gamma<n$, $1\leq p<n/\gamma$ and $1/q=1/p-\gamma/n$. Let $\Phi,\xi$ be Young functions verifying $t^{\gamma/n}\xi^{-1}(t)\leq C \Phi^{-1}(t)$, for every $t\geq t_0\geq 0$. Then, for every nonnegative functions $w$ and $f_0\in L^p$ we have that
	\[M_{\gamma,\Phi}\left(\frac{f_0}{w}\right)(x)\leq CM_\xi\left(\frac{f_0^{p/q}}{w}\right)(x)\left(\int_{\mathbb{R}^n}f_0^p(y)\,dy\right)^{\gamma/n},\]
	for every $x\in \mathbb{R}^n$. 
\end{propo}

\begin{proof}[Proof of Theorem~\ref{teo: mixta para Mgamma,phi}]
	Assume without loss of generality that the right-hand side is finite. We take $\nu=\delta q/r$ and define $\xi(t)=t^q(1+\log^+t)^\nu$. Notice that by \eqref{eq: preliminares - inversa de L log L} we have that
	\[t^{\gamma/n}\xi^{-1}(t)\approx \frac{t^{\gamma/n+1/q}}{(1+\log^+t)^{\delta/r}}\approx \Phi^{-1}(t).\]
	 If we set $f_0=|f|v(Mu)^{1/q}$ we have that $f_0\in L^r$ by the assumption above. Therefore we apply Proposition~\ref{propo: estimacion puntual M_{gamma, Phi}} with $p=r$ and $w=(Mu)^{1/q}$ to get
	\[M_{\gamma,\Phi}\left(\frac{f_0}{w}\right)(x)\leq CM_\xi\left(\frac{f_0^{r/q}}{w}\right)(x)\left(\int_{\mathbb{R}^n}f_0^r(y)\,dy\right)^{\gamma/n}.\]
	Hence
	\begin{align*}
	uv^{q}\left(\left\{x: \frac{M_{\gamma,\Phi}(fv)(x)}{v(x)}>t\right\}\right)&=uv^{q}\left(\left\{x: \frac{M_{\gamma,\Phi}(f_0/w)(x)}{v(x)}>t\right\}\right)\\
	&\leq uv^{q}\left(\left\{x: \frac{M_{\xi}(f_0^{r/q}/w)(x)}{ v(x)}>\frac{t}{\left(\int f_0^r \right)^{\gamma/n}}\right\}\right).
	\end{align*}
	Since $\xi\in \mathfrak{F}_q$,  we can use the mixed estimate for $M_\xi$ given by Theorem~\ref{teo: desigualdad mixta para MPhi} to obtain
	\begin{equation}\label{eq: eq1 - teo mixta para M_{gamma,Phi}, caso p=r}
	uv^{q}\left(\left\{x: \frac{M_{\gamma,\Phi}(fv)(x)}{v(x)}>t\right\}\right)\leq C\int_{\mathbb{R}^n}\xi\left(\frac{f_0^{r/q}\left(\int f_0^r\right)^{\gamma/n}}{wv t}\right)v^qMu.
	\end{equation}
	Recalling the relation $1/q=1/r-\gamma/n$ or equivalently $r/q+r\gamma/n=1$, the argument of $\xi$ can be written as
	
	\begin{align*}
	\frac{f_0^{r/q}\left(\int f_0^r\right)^{\gamma/n}}{wv t}&=\left(\frac{|f|}{t}\right)^{r/q}(wv)^{r/q-1}\left(\int_{\mathbb{R}^n} \left(\frac{|f|}{t}\right)^r(wv)^r\right)^{\gamma/n}\\
	&=\left[\left(\frac{|f|}{t}\right)(wv)^{1-q/r}\left(\int_{\mathbb{R}^n} \left(\frac{|f|}{t}\right)^r(wv)^r\right)^{\gamma q/(nr)}\right]^{r/q}.
	\end{align*}
	Notice that $\xi(t^{r/q})=t^r$ for $0\leq t\leq 1$. On the other hand, if $t>1$ we have
	\begin{align*}
	\xi(t^{r/q})&=t^r(1+\log t^{r/q})^\nu\\
	&=t^r\left(1+\frac{r}{q}\log t\right)^\nu,
	\end{align*}
	which implies $\xi(t^{r/q})\leq \Phi_\gamma(t)=t^r(1+\log^+t)^\nu$. Since $\Phi_\gamma$ is submultiplicative, we proceed now as follows
	\begin{align*}
	\xi\left(\frac{f_0^{r/q}\left(\int_{\mathbb{R}^n}f_0^r\right)^{\gamma/n}}{wv t}\right)&\leq \Phi_\gamma\left(\left(\frac{|f|}{t}\right)(wv)^{1-q/r}\left(\int_{\mathbb{R}^n} \left(\frac{|f|}{t}\right)^r(wv)^r\right)^{\gamma q/(nr)}\right)\\
	&\leq \Phi_\gamma\left(\left[\int_{\mathbb{R}^n}\Phi_\gamma\left(\frac{|f|}{t}\right)(wv)^r\right]^{\gamma q/(nr)}\right)\Phi_\gamma\left(\frac{|f|}{t}(wv)^{1-q/r}\right).
	\end{align*}
	Returning to \eqref{eq: eq1 - teo mixta para M_{gamma,Phi}, caso p=r}, the right hand side is bounded by
	\[ \Phi_\gamma\left(\left[\int_{\mathbb{R}^n}\Phi_\gamma\left(\frac{|f|}{t}\right)(wv)^r\right]^{\gamma q/(nr)}\right)\int_{\mathbb{R}^n} \Phi_\gamma\left(\frac{|f|}{t}(wv)^{1-q/r}\right)(wv)^q.\]
	Observe that $\Phi_\gamma(t^{1-q/r})t^q\leq \Psi(t)=t^r(1+\log^+(t^{1-q/r}))^{n\delta/(n-r\gamma)}$. Therefore, the expression above is bounded by 
	\[\Phi_\gamma\left(\left[\int_{\mathbb{R}^n}\Phi_\gamma\left(\frac{|f|}{t}\right)\Psi((Mu)^{1/q}v)\right]^{\gamma q/(nr)}\right)\int_{\mathbb{R}^n} \Phi_\gamma\left(\frac{|f|}{t}\right)\Psi((Mu)^{1/q}v).\]
	In order to conclude, notice that 
	\[t\Phi_\gamma(t^{\gamma q/(nr)})\lesssim t^{1+\gamma q/n}(1+\log^+ t)^\nu= t^{q/r}(1+\log^+ t)^{\delta q/r}=\varphi(t).\qedhere\]
\end{proof}

\section{A particular case of special interest}\label{seccion: caso Phi(t)=t}

In this section we give the proof of Theorem~\ref{teo: mixta para I_gamma}. We shall need the following extrapolation result. 

\begin{teo}[\cite{O-P}, Theorem 1.1]\label{teo: extrapolacion}
	Let $\mathcal{F}$ be a family of pairs of nonnegative functions $(f,g)$ that verify the inequality
	\begin{equation}\label{eq: tipo coifman entre f y g}
	 \int_{\mathbb{R}^n} f(x)^{p_0}w(x)\,dx\leq C\int_{\mathbb{R}^n}g(x)^{p_0}w(x)\,dx
	\end{equation}
	 for some $0<p_0<\infty$ and every $w\in A_\infty$. If $\theta\geq 1$, $u\in A_1$ and $v$ is a function such that $v^\delta\in A_\infty$ for some $\delta>0$, then we have that
	 \[\left\|\frac{f}{v^\theta}\right\|_{L^{1/\theta,\infty}(uv)}\leq C\left\|\frac{g}{v^\theta}\right\|_{L^{1/\theta,\infty}(uv)}.\] 
\end{teo}

The key estimate to prove mixed inequalities for $I_\gamma$ is a consequence of the theorem above. 

\begin{propo}\label{propo: control de I_gamma por M_gamma}
	Let $0<\gamma<n$, $q=n/(n-\gamma)$, $u\in A_1$ and $v(x)=|x|^\beta$ with $\beta<-n$. Then there exists a positive constant $C$ such that
	 \[\left\|\frac{I_\gamma(fv)}{v}\right\|_{L^{q,\infty}(uv^q)}\leq C\left\|\frac{M_\gamma(fv)}{v}\right\|_{L^{q,\infty}(uv^q)}.\]	  
\end{propo}

\begin{proof}
	Let $F=|I_\gamma(fv)|^q$ and $G=(M_\gamma(fv))^q$. By using \eqref{eq: intro - tipo Coifman Igamma-Mgamma} we have that \eqref{eq: tipo coifman entre f y g} holds. Notice that $|x|^\alpha$ belongs to $A_1$ if and only if $-n<\alpha<0$. If we pick $0<\delta<-n/(q\beta)$ we have that $v^{q\delta}\in A_1\subseteq A_\infty$. By applying Theorem~\ref{teo: extrapolacion} with $u$, $v^q$ and $\theta=1$ we get
	\[\left\|\frac{F}{v^q}\right\|_{L^{1,\infty}(uv^q)}\leq C\left\|\frac{G}{v^q}\right\|_{L^{1,\infty}(uv^q)}.\]
	Hence we can proceed as follows
	\begin{align*}
	\left\|I_\gamma(fv)v^{-1}\right\|_{L^{q,\infty}(uv^q)}&=\left\|\left(Fv^{-q}\right)^{1/q}\right\|_{L^{q,\infty}(uv^q)}\\
	&=\left\|Fv^{-q}\right\|_{L^{1,\infty}(uv^q)}^{1/q}\\
	&\leq C\left\|Gv^{-q}\right\|_{L^{1,\infty}(uv^q)}^{1/q}\\
	&=\left\|M_\gamma(fv)v^{-1}\right\|_{L^{q,\infty}(uv^q)}.\qedhere
	\end{align*}
\end{proof}

\begin{proof}[Proof of Theorem~\ref{teo: mixta para I_gamma}]
	The estimate follows by combining Corollary~\ref{coro: mixta para M_gamma} with Proposition~\ref{propo: control de I_gamma por M_gamma}. Indeed, for every $t>0$ we have that
	\begin{align*}
	tuv^q\left(\left\{x\in\mathbb{R}^n: \frac{|I_\gamma(fv)(x)|}{v(x)}>t\right\}\right)^{1/q}&\leq \sup_{t>0}t\,uv^q\left(\left\{x\in\mathbb{R}^n: \frac{|I_\gamma(fv)(x)|}{v(x)}>t\right\}\right)^{1/q}\\
	&=\left\|\frac{I_\gamma(fv)}{v}\right\|_{L^{q,\infty}(uv^q)}\\
	&\leq C\left\|\frac{M_\gamma(fv)}{v}\right\|_{L^{q,\infty}(uv^q)}\\
	&=C\sup_{t>0}\left[t\,uv^q\left(\left\{x\in\mathbb{R}^n: \frac{M_\gamma(fv)(x)}{v(x)}>t\right\}\right)^{1/q}\right]\\
	&\leq C\int_{\mathbb{R}^n}|f(x)|u^{1/q}(x)v(x)\,dx,
	\end{align*}
	which implies the thesis.
\end{proof}

\section{Further applications: commutators with Lipschitz symbol}\label{seccion: aplicaciones}

As a final application we can give mixed inequalities for commutators of certain operators with Lipschitz symbol. We shall be dealing with a linear operator $T$ bounded on $L^2(\R^n)$ and such that for $f\in L^2$ with compact support we have the representation
\begin{equation}\label{eq: representacion integral del T}
Tf(x)=\int_{\R^n}K(x,y)f(y)\,dy ,\quad\quad x\notin \supp f,
\end{equation}
where $K$ is a measurable function defined away from the origin.
Recall that $T$ is a \emph{Calderón-Zygmund operator} (CZO) if $K$ is a standard kernel, which means that
$K:\mathbb R^n\times \mathbb R^n\backslash\Delta\to\mathbb C$  satisfies a size condition given by 
\[|K(x,y)|\lesssim \frac{1}{|x-y|^n},\]
and the smoothness conditions, usually called Lipschitz conditions,
\begin{align*}
|K(x,y)-K(x,z)|&\lesssim \frac{|x-z|}{|x-y|^{n+1}},\quad \textrm{ if } |x-y|>2|y-z|\\
|K(x,y)-K(w,z)|&\lesssim\frac{|x-w|}{|x-y|^{n+1}},\quad \textrm{ if } |x-y|>2|x-w|.
\end{align*}

We say that $T$ is a \emph{Calderón-Zygmund singular integral} operator if
\[Tf(x)=\lim_{\varepsilon\to 0} \int_{|y-x|>\epsilon}K(x,y)f(y)\,dy=\lim_{\varepsilon\to 0} T_\varepsilon f(x).\]

Given an operator $T$ and a locally integrable function $b$ the commutator of $T$ is denoted by $T_b$ or $[b,T]$ and defined by the expression
\[[b,T]f=bTf-T(bf).\]
For $m\in \mathbb{N}$, the higher order commutator of $T$ is given recursively by $T_b^m=[b, T_b^{m-1}]$.

When $T=I_\gamma$ we will denote the commutator of order $m$ by $I_{\gamma,b}^m$.

 Given a locally integrable function $b$ and $0<\delta\leq 1$ we say that $b$ belongs to the class $\Lambda(\delta)$ if there exists a positive constant $C$ such that
\[|b(x)-b(y)|\leq C|x-y|^\delta\]
for every $x,y\in \mathbb{R}^n$. The smallest constant $C$ for which the expression above holds will be denoted by $\|b\|_{\Lambda(\delta)}$.

The following two lemmas contain a pointwise relation between commutators with Lipschitz symbol and the fractional integral operator. Both were set and proved in \cite{BCP-J}. 

\begin{lema}\label{lema: relacion puntual para conmutadores de CZSIO}
	Let $m\in\mathbb{N}$, $0<\delta<\min\{1, n/m\}$, $b\in\Lambda(\delta)$ and $T$ a Calderón-Zygmund singular integral operator. For almost every $x\in\mathbb{R}^n$ and every bounded function $f$ with compact support we have that
	\[|T_b^mf(x)|\leq C\|b\|_{\Lambda(\delta)}^mI_{m\delta}|f|(x).\]
\end{lema}

\begin{lema}\label{lema: relacion puntual para conmutadores de I_gamma}
	Let $0<\gamma<n$, $m\in\mathbb{N}$, $0<\delta<\min\{1, (n-\gamma)/m\}$ and $b\in\Lambda(\delta)$. Then the inequality
	\[|I_{\gamma,b}^mf(x)|\leq C\|b\|_{\Lambda(\delta)}^mI_{\gamma+m\delta}|f|(x)\]
	holds for almost every $x$ and every bounded function $f$ with compact support.
\end{lema}

By combining these estimates with Theorem~\ref{teo: mixta para I_gamma} we can obtain the two following results.

\begin{teo}
	Let $m\in\mathbb{N}$,  $0<\delta<\min\{1, n/m\}$, $b\in\Lambda(\delta)$ and $T$ a Calderón-Zygmund singular integral operator. Let $q=n/(n-m\delta)$, $u\in A_1$ and $v(x)=|x|^\beta$ with $\beta<-n$. Then there exists a positive constant $C$ such that the inequality
\[uv^q\left(\left\{x\in\mathbb{R}^n: \frac{|T_b^m(fv)(x)|}{v(x)}\right\}\right)\leq \frac{C}{t}\|b\|_{\Lambda(\delta)}^m\int_{\mathbb{R}^n}|f(x)|u^{1/q}(x)v(x)\,dx\]
holds for every $t>0$ and every bounded function $f$ with compact support.
\end{teo}

\begin{teo}
	Let $m\in\mathbb{N}$, $0<\gamma<n$, $0<\delta<\min\{1,(n-\gamma)/m\}$ and $b\in\Lambda(\delta)$. If $q=n/(n-\gamma-m\delta)$, $u\in A_1$ and $v(x)=|x|^\beta$ with $\beta<-n$, there exists a positive constant $C$ such that the inequality
	\[uv^q\left(\left\{x\in\mathbb{R}^n: \frac{|I_{\gamma,b}^m(fv)(x)|}{v(x)}\right\}\right)\leq \frac{C}{t}\|b\|_{\Lambda(\delta)}^m\int_{\mathbb{R}^n}|f(x)|u^{1/q}(x)v(x)\,dx\]
	holds for every $t>0$ and every bounded function $f$ with compact support.
\end{teo}

\section*{Acknowledgements}
The author was supported by CONICET and UNL.


\def\cprime{$'$}
\providecommand{\bysame}{\leavevmode\hbox to3em{\hrulefill}\thinspace}
\providecommand{\MR}{\relax\ifhmode\unskip\space\fi MR }
\providecommand{\MRhref}[2]{%
	\href{http://www.ams.org/mathscinet-getitem?mr=#1}{#2}
}
\providecommand{\href}[2]{#2}

\end{document}